\documentclass[12pt,reqno]{amsart}
\def\draft{n}
\usepackage{
  graphicx, amsmath, amssymb, multicol, mathtools,
  picins, stmaryrd, xcolor, mathabx, needspace, enumitem,
  lgreek, sidecap, mathbbol, blindtext, pdfpages
}
\usepackage[all]{xy}
\usepackage[textwidth=6.5in,textheight=9in,headsep=0.15in,centering]{geometry}

\usepackage[breaklinks=true]{hyperref}\hypersetup{colorlinks,
  linkcolor={green!50!black},
  citecolor={green!50!black},
  urlcolor=blue
}

\usepackage{chngcntr}
\counterwithin{figure}{section}

\usepackage[within=section]{newfloat}
\DeclareFloatingEnvironment[placement=b]{Aside}
\usepackage{caption}
\if\draft y
  \usepackage[color]{showkeys}
\fi

\theoremstyle{plain}
\newtheorem{theorem}{Theorem}

\theoremstyle{definition}

\newtheorem{note}[theorem]{Note}

\newlength{\standardunitlength}
\def\qed{{\linebreak[1]\null\hfill\text{$\Box$}}}

\newlength{\globalparindent}
\def\arXiv#1{{\href{https://arxiv.org/abs/#1}{arXiv:\linebreak[0]#1}}}

\def\weburl{http://drorbn.net/SL2PO}

\def\bbR{{\mathbb R}}

\def\fraka{{\mathfrak a}}

\def\frakg{{\mathfrak g}}
\def\frakh{{\mathfrak h}}
\def\frakl{{\mathfrak l}}

\def\fraku{{\mathfrak u}}

\def\draftcut{\if\draft y \cleardoublepage \fi}

\definecolor{lightred}{RGB}{255, 217, 217}

\def\eps{\epsilon}

\begin{document} 
\newdimen\captionwidth\captionwidth=\hsize
\setcounter{secnumdepth}{4}

\title{An Unexpected Cyclic Symmetry of $I\fraku_n$}

\author{Dror~Bar-Natan}
\address{
  Department of Mathematics\\
  University of Toronto\\
  Toronto Ontario M5S 2E4\\
  Canada
}
\email{drorbn@math.toronto.edu}
\urladdr{http://www.math.toronto.edu/drorbn}

\author{Roland~van~der~Veen}
\address{
  University of Groningen, Bernoulli Institute\\
  P.O. Box 407\\
  9700 AK Groningen\\
  The Netherlands
}
\email{roland.mathematics@gmail.com}
\urladdr{http://www.rolandvdv.nl/}

\date{First edition February 3, 2020, this edition \today.}

\subjclass[2010]{57M25}
\keywords{
  Lie algebras,
  Lie bialgebras,
  Lie algebra automorphism,
  solvable approximation,
  triangular matrices}

\thanks{This work was partially supported by NSERC grant RGPIN-2018-04350. It is available in electronic form, along with source files and a verification {\sl Mathematica} notebook at \url{http://drorbn.net/UnexpectedCyclic} and at \arXiv{2002.00697}.}

\begin{abstract}
We find and discuss an unexpected (to us) order $n$ cyclic group of automorphisms of the Lie algebra $I\fraku_n\coloneqq\fraku_n\ltimes\fraku_n^\ast$, where $\fraku_n$ is the Lie algebra of upper triangular $n\times n$ matrices. Our results also extend to $gl_{n+}^\eps$, a ``solvable approximation'' of $gl_n$, as defined within.
\end{abstract}

\maketitle


Given any Lie algebra $\fraka$ one may form its ``inhomogeneous version'' $I\fraka\coloneqq\fraka\ltimes\fraka^\ast$, its semidirect product with its dual $\fraka^\ast$ where $\fraka^\ast$ is considered as an Abelian Lie\footnote{Two Norwegians!} algebra and $\fraka$ acts on $\fraka^\ast$ via the coadjoint action. (Over $\bbR$ if $\fraka=so_3$ then $\fraka^\ast=\bbR^3$ and so $I\fraka=so_3\ltimes\bbR^3$ is the Lie algebra of the Euclidean group of rotations and translations, explaining the name).

In general, we care about $I\fraka$. It is a special case of the Drinfel'd double / Manin triple construction~\cite{Drinfeld:QuantumGroups, EtingofSchiffman:QuantumGroups} when the cobracket is $0$. These Lie algebras occur in the study of the Kashiwara-Vergne problem~\cite{Talk:Bonn, WKO2} and they provide the simplest quantum algebra context for the Alexander polynomial~\cite{Talk:Chicago, WKO1}. We care especially for the case where $\fraka$ is a Borel subalgebra of a semi-simple Lie algebra (e.g., upper triangular matrices) as then the algebras $I\fraka$ are the $\eps=0$ ``base case'' for ``solvable approximation''~\cite{PP1, DoPeGDO, DPG, Talk:SolvApp, Talk:Dogma, Talk:DoPeGDO}, and their automorphisms are expected to become symmetries of the resulting knot invariants.

\begin{figure}
\[ \resizebox{!}{3cm}{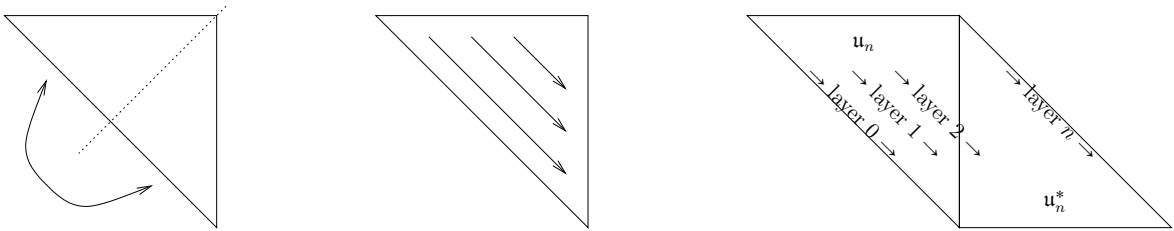} \]
\caption{An expected anti-automorphism (left), an unexpected automorphism (middle), and an alternative presentation of the ``layers'' table (right).} \label{fig:ExpectedAndUnexpected}
\end{figure}

Let $\fraku_n$ be the Lie algebra of upper triangular $n\times n$ matrices. Beyond inner automorphisms, $\fraku_n$ and hence $I\fraku_n$ has one obvious and expected anti-automorphism $\Phi$ corresponding to flipping matrices along their anti-main-diagonal, as shown in the first image of Figure~\ref{fig:ExpectedAndUnexpected}. With $x_{ij}$ denoting the $n\times n$ matrix with $1$ in position $(ij)$ and zero everywhere else ($i\leq j$ in $\fraku_n$), $\Phi$ is given by $x_{ij}\mapsto x_{n+1-j,n+1-i}$.

There clearly isn't an automorphism of $\fraku_n$ that acts by ``sliding down and right parallel to the main diagonal'', as in the second image in Figure~\ref{fig:ExpectedAndUnexpected}. Where would the last column go? Yet the sliding map, when restricted to where it is clearly defined ($\fraku_n$ with the last column excluded), does extend to an automorphism of $I\fraku_n$ as in the theorem below.

\begin{theorem} \label{thm:Psi4Iu} With the basis $\{x_{ij}\}_{1\leq i<j\leq n} \cup \{a_i=x_{ii}\}_{1\leq i\leq n}$ for $\fraku_n$ and dual basis $\{x_{ji}\}_{1\leq i<j\leq n}\cup\{b_i\}_{1\leq i\leq n}$ for $\fraku_n^\ast$ (and duality $\langle x_{kl},x_{ij}\rangle=\delta_{li}\delta_{jk}$, $\langle b_i,a_j\rangle=2\delta_{ij}$,\footnotemark\ and $\langle x_{ji},a_k\rangle = \langle b_k,x_{ij}\rangle = 0$), the map $\Psi\colon I\fraku_n\to I\fraku_n$ defined by ``incrementing all indices by $1$ mod $n$'' (precisely, if $\psi$ is the single-cycle permutation $\psi=(123\ldots n)$ then $\Psi$ is defined by $\Psi(x_{ij})=x_{\psi(i)\psi(j)}$, $\Psi(a_i)=a_{\psi(i)}$, and $\Psi(b_i)=b_{\psi(i)}$) is a Lie algebra automorphism of $I\fraku_n$.
\end{theorem}

\footnotetext{The awkward factor of $2$ in $\langle b_i,a_j\rangle$ is irrelevant for Theorem~\ref{thm:Psi4Iu} yet crucial for Theorem~\ref{thm:Psi4gleps}. Removing this factor removes the factor $\frac12$ in~\eqref{eq:cr0}.}

Note that our choice of bases, using similar symbols $x_{ij}$ / $x_{ji}$ for the non-diagonal matrices and their duals, hides the intricacy of $\Psi$; e.g., $\Psi\colon x_{n-1,n}\mapsto x_{n1}$ maps an element of $\fraku_n$ to an element of $\fraku_n^\ast$ (also see Figure~\ref{fig:ExpectedAndUnexpected}, right).

It may be tempting to think that $\Psi$ has a simple explanation in $gl_n$ language: $\fraku_n$ is a subset of $gl_n$, $gl_n$ has a metric (the Killing form) such that the dual of $x_{ij}$ is $x_{ji}$ as is the case for us, and every permutation of the indices induces an automorphism of $gl_n$. But this explains nothing and too much: nothing because the bracket of $I\fraku_n$ simply isn't the bracket of $gl_n$ (even away from the diagonal matrices), and too much because {\em every} permutation of indices induces an automorphism of $gl_n$, whereas only $\psi$ and its powers induce automorphisms of $I\fraku_n$.

\vskip 1mm
\noindent{\bf Proof of Theorem~\ref{thm:Psi4Iu}.} Recall that as a vector space $I\fraku_n=\fraku_n\oplus\fraku_n^\ast$, yet with bracket $[(x,f),(y,g)]=\left([x,y],x\cdot g-y\cdot f\right)$ where $\cdot$ denotes the coadjoint action, $(x\cdot f)(v)=f([v,x])$. With that and some case checking and explicit computations, the commutation relations of $I\fraku_n$ are given by
\begin{equation} \label{eq:cr0}
  \begin{aligned}
    [x_{ij},x_{kl}] &= \chi_{\lambda(x_{ij})+\lambda(x_{kl})<n}(\delta_{jk}x_{il}-\delta_{li}x_{kj})
      & \text{unless both }j=k\text{ and }l=i, \\
    [x_{ij},x_{ji}] &= \frac12(b_i-b_j), & \\
    [a_i,x_{jk}] &= (\delta_{ij}-\delta_{ik})x_{jk}, & \\
    [b_i,x_{jk}] &= 0, & \\
    [a_i,a_j] &= [b_i,b_j] = [a_i,b_j] = 0, &
  \end{aligned}
\end{equation}
where $\chi$ is the indicator function of truth, $\chi_{5<7}=1$ while $\chi_{7<5}=0$, and where $\lambda(x_{ij})$ is the ``length'' of $x_{ij}$, defined by $\lambda(x_{ij})\coloneqq\begin{cases} j-i & i<j \\ n-(i-j) & i>j \end{cases}$.

It is easy to verify that the length $\lambda(x_{ij})$ is $\Psi$-invariant, and hence everything in~\eqref{eq:cr0} is $\Psi$-equivariant. \qed

$I\fraku_n$ is a solvable Lie algebra (as a semi-direct product of solvable with Abelian, and as will be obvious from the table below). It is therefore interesting to look at the structure of its commutator subgroups. This structure is summarized in the following table (an alternative view is in Figure~\ref{fig:ExpectedAndUnexpected}):

\needspace{45mm}
\[ \arraycolsep=2pt \def\arraystretch{1.2}
  \begin{array}{|l|c|cccccccccc|}
    \hline
    \text{layer }0 & \frakg=I\fraku_n & a_1 & \to & a_2  & \to \cdots \to & a_{n-2} & \to & a_{n-1} & \to & a_n & \to \\
    \text{layer }1 & \frakg'_1=\frakg'=[\frakg,\frakg] & x_{12} & \to & x_{23} & \to \cdots \to & x_{n-2,n-1} & \to & x_{n-1,n} & \to & x_{n1} & \to \\
    \text{layer }2 & \frakg'_2=[\frakg',\frakg'_1] & x_{13} & \to & x_{24} & \to \cdots \to & x_{n-2,n} & \to & x_{n-1,1} & \to & x_{n2} & \to \\
    \text{layer }3 & \frakg'_3=[\frakg',\frakg'_2] & x_{14} & \to & x_{25} & \to \cdots \to & x_{n-2,1} & \to & x_{n-1,2} & \to & x_{n3} & \to \\
    \vdots & \vdots & \vdots & & \vdots & & \vdots & & \vdots & & \vdots & \\
    \text{layer }(n-1) & \frakg'_{n-1}=[\frakg',\frakg'_{n-2}] & x_{1n} & \to & x_{21} & \to \cdots \to & x_{n-2,n-1} & \to & x_{n-1,n-2} & \to & x_{n,n-1} & \to \\
    \text{layer }n & \frakg'_n=[\frakg',\frakg'_{n-1}] & b_1 & \to & b_2 & \to \cdots \to & b_{n-2} & \to & b_{n-1} & \to & b_n & \to \\
    \hline
  \end{array}
\]
In this table (all assertions are easy to verify):
\begin{itemize}[leftmargin=*,labelindent=0pt]
\item Apart for the treatment of the $a_i$'s and the $b_i$'s, layer$=$length$=\lambda(x_{ij})$.
\item The layers indicate a filtration; each layer should be considered to contain all the ones below it. The generators marked at each layer generate it modulo the layers below.
\item The bracket of an element at layer $p$ with an element of layer $q$ is in layer $p+q$ (and it must vanish if $p+q>n$).
\item If $p\geq 2$, every generator in layer $p$ is the bracket of a generator in layer $1$ with a generator in layer $p-1$.
\item In layer $p$, the first $n-p$ generators indicated belong to $\fraku_n$ and the last $p$ belong to $\fraku_n^\ast$. So as we go down, $\fraku_n^\ast$ slowly ``overtakes'' the table.
\item The automorphism $\Psi$ acts by following the arrows and shifting every generator one step to the right (and pushing the rightmost generator in each layer back to the left).
\item The anti-automorphism $\Phi$ acts by mirroring the $\fraku_n$ part of each layer left to right and by doing the same to the $\fraku_n^\ast$ part, without mixing the two parts.
\item Note that $I\fraku_n$ can be metrized by pairing the $\fraku_n$ summand with the $\fraku_n^\ast$ one. The metric only pairs generators indicated in layer $p$ with generators indicated in layer $(n-p)$.
\end{itemize}

Note also that the brackets of the generators indicated in layer 1 yield the generators indicated in layer 2 as follows:
\[ \xymatrix@R=10pt@C=10pt{
  x_{12} \ar[rd] && x_{23} \ar[ld]\ar[rd] && x_{34} \ar[ld]\ar[rd] && \cdots && x_{n-1,n} \ar[ld]\ar[rd] && x_{n1} \ar[ld]\ar[rd] && \ar[ld] \\
  & x_{13} && x_{24} &&& \cdots &&& x_{n-1,1} && x_{n2} &
} \]
(with the diagram continued cyclically). The symmetry group of the above cycle is the dihedral group $D_n$ and this strongly suggests that the group of outer automorphisms and anti-automorphisms of $I\fraku_n$ (all automorphisms and anti-automorphisms modulo inner automorphisms) is $D_n$, generated by $\Phi$ and $\Psi$. We did not endeavor to prove this formally.

\vskip 2mm\noindent{\bf Extension.} The Drinfel'd double / Manin triple construction~\cite{Drinfeld:QuantumGroups, EtingofSchiffman:QuantumGroups}, when applied to $\fraku_n$, is a way to reconstruct $gl_n$ from its subalgebras of upper triangular matrices $\fraku_n$ and lower triangular matrices $\frakl_n$. Precisely, one endows the vector space $\frakg=\fraku_n\oplus\frakl_n$ with a non-degenerate symmetric bilinear form by declaring that the subspaces $\fraku_n$ and $\frakl_n$ are isotropic ($\langle\fraku_n,\fraku_n\rangle=\langle\frakl_n,\frakl_n\rangle=0$) and by setting $\langle x_{kl},x_{ij}\rangle=\delta_{li}\delta_{jk}$, $\langle b_i,a_j\rangle=2\delta_{ij}$, and $\langle x_{ji},a_k\rangle = \langle b_k,x_{ij}\rangle = 0$ as in Theorem~\ref{thm:Psi4Iu} and where $a_i$ stands for the diagonal matrix $x_{ii}$ considered as an element of $\fraku_n$ and $b_i$ stands for the same matrix as an element of $\frakl_n$. There is then a unique bracket on $\frakg$ that extends the brackets on the summands $\fraku_n$ and $\frakl_n$ and relative to which the inner product of $\frakg$ is invariant. With our judicious choice of bilinear form, this bracket on $\frakg$ satisfies the Jacobi identity and turns $\frakg$ into a Lie algebra isomorphic to $gl_{n+}=gl_n\oplus\frakh'_n$, where $\frakh'_n$ denotes a second copy of the diagonal matrices in $gl_n$.

We let $gl_{n+}^\eps$ be the Inonu-Wigner~\cite{InonuWigner:ContractionOfGroups} contraction of $\frakg$ along its $\frakl_n$ summand, with parameter $\epsilon$.\footnote{Alternatively, make $\fraku_n$ into a Lie bialgebra with cobracket $\delta$ using its given duality with $\frakl_n$, and double it as in~\cite{Drinfeld:QuantumGroups, EtingofSchiffman:QuantumGroups} but using the cobracket $\eps\delta$.} All that this means is that the bracket of $\frakl_n$ gets multiplied by $\epsilon$ to give $\frakl_n^\epsilon$, and then the Drinfel'd double / Manin triple construction is repeated starting with $\fraku_n\oplus\frakl_n^\eps$, without changing the bilinear form. The result is a Lie algebra $gl_{n+}^\eps$ over the ring of polynomials in $\eps$ which specializes to $I\fraku_n$ at $\eps=0$ and which is isomorphic to $gl_n\oplus\frakh'_n$ when $\eps$ is invertible.\footnote{Hence $gl_{n+}^\eps\to I\fraku_n$ is a counter-example to the feel-true statement ``a contraction of a direct sum is a direct sum''. Indeed with notation as in Theorem~\ref{thm:Psi4gleps}, as $\eps\to 0$ the decomposition $gl_{n+}^\eps = gl_n\oplus\frakh'_n = \langle x_{ij},b_i+\eps a_i\rangle \oplus \langle b_i-\eps a_i\rangle$ collapses.}\ We care about $gl_{n+}^\eps$ a lot~\cite{PP1, DoPeGDO, DPG, Talk:SolvApp, Talk:Dogma, Talk:DoPeGDO}; when reduced modulo $\eps^{k+1}=0$ for some natural number $k$ it becomes solvable, and hence a ``solvable approximation'' of $gl_n$ with applications to computability of knot invariants.

\begin{theorem} \label{thm:Psi4gleps} With the same conventions as in Theorem~\ref{thm:Psi4Iu} the map $\Psi$ is also a Lie algebra automorphism of $gl_{n+}^\eps$.
\end{theorem}

\begin{proof} By some case checking and explicit computations, the commutation relations of $gl_{n+}^\eps$ are given by
\[
  \begin{aligned}
    [x_{ij},x_{kl}] &= \chi^\eps_{\lambda(x_{ij})+\lambda(x_{kl})<n}(\delta_{jk}x_{il}-\delta_{li}x_{kj})
      & \text{unless both }j=k\text{ and }l=i, \\
    [x_{ij},x_{ji}] &= \frac12(b_i-b_j)+\frac{\eps}{2}(a_i-a_j), & \\
    [a_i,x_{jk}] &= (\delta_{ij}-\delta_{ik})x_{jk}, & \\
    [b_i,x_{jk}] &= \eps(\delta_{ij}-\delta_{ik})x_{jk}, & \\
    [a_i,a_j] &= [b_i,b_j] = [a_i,b_j] = 0, &
  \end{aligned}
\]
where $\chi^\eps_{\scriptsize\verb"True"}=1$ and $\chi^\eps_{\scriptsize\verb"False"}=\eps$. These relations are clearly $\Psi$-equivariant. \qed
\end{proof}

\begin{note} There is of course an ``$sl$'' version of everything, in which linear combinations $\sum\alpha_ia_i$ and $\sum\beta_ib_i$ are allowed only if $\sum\alpha_i=\sum\beta_i=0$, with obvious modifications throughout.
\end{note}

\begin{note} At $n=2$ and $\eps=0$, the algebra $sl_{2+}^0$ is the ``diamond Lie algebra'' of~\cite[Chapter~4.3]{Kirillov:OrbitMethod}, which is sometimes called ``the Nappi-Witten algebra''~\cite{NappiWitten:Nonsemisimple}: With $a=(a_1-a_2)/2$, $x=x_{12}$, $y=x_{21}$, and $b=(b_1-b_2)/2$, it is
\[ \langle a,x,y,b\rangle/\left([a,x]=x,\,[a,y]=-y,\,[x,y]=b,\,[b,-]=0\right). \]
Here $\Phi\colon(a,x,y,b)\mapsto(-a,x,y,-b)$ and $\Psi\colon(a,x,y,b)\mapsto(-a,y,x,-b)$.
\end{note}

\begin{note} \label{note:KZ} Upon circulating this paper as an eprint we received a note from A.~Knutson informing us of~\cite[esp.~sec.~2.3]{KnutsonZinnJustin:BrauerLoopModel}, where the algebra $I\fraku_n$ (except reduced modulo $\langle b_i\rangle$ and considered globally rather than infinitesimally) is considered from a different perspective. It is shown to be a subquotient of the affine algebra $\widehat{gl_n}$ in a manner preserved by its automorphisms corresponding to its Dynkin diagram, which is a cycle. Similar comments apply to the other algebras considered here.
\end{note}

\begin{note} A day later we received a note~\cite{BuloisRessayre:Automorphisms} from M.~Bulois and N.~Ressayre reporting on an explanation of Theorem~\ref{thm:Psi4Iu} in terms of affine Kac-Moody Lie algebras, similarly to Note~\ref{note:KZ}.
\end{note}

\vskip 2mm\noindent{\bf Acknowledgement.} We wish to thank M.~Bulois, A.~Knutsen, N.~Ressayre, and N.~Williams for comments and suggestions.

\parpic[l]{\includegraphics[height=2in]{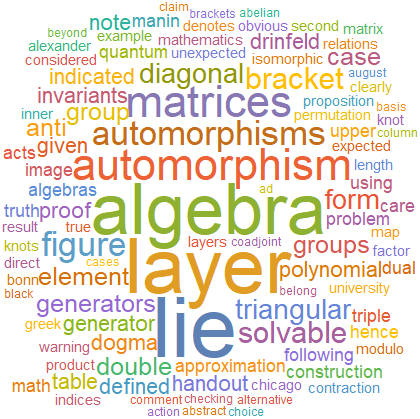}}

\AtEndDocument{

\section*{Supplementary material}

# An Unexpected Cyclic Symmetry of $Iu_n$ - Verification Notebook

Pensieve header: Verification notebook for "An Unexpected Cyclic Symmetry of $Iu_n$ by Dror Bar-Natan and Roland van der Veen. Also available at http://drorbn.net/UnexpectedCyclic. Continues pensieve://2020-01/ .

Only Theorem 2 is tested; Theorem 1 is simply the case where $\epsilon = 0$, so it does not require independent testing.

## Definitions.

General definitions - brackets $B$ and pairings $P$ are bilinear, brackets are anti-symmetric:

```
In[◦]:= B[0, _] = 0; B[_, 0] = 0;
       B[c_ * x : (x | a | b)__, y_] := Expand[c B[x, y]];
       B[y_, c_ * x : (x | a | b)__] := Expand[c B[y, x]];
       B[x_Plus, y_] := B[#, y] & /@ x;
       B[x_, y_Plus] := B[x, #] & /@ y;

In[◦]:= P[0, _] = 0; P[_, 0] = 0;
       P[c_ * x : (x | a | b)__, y_] := Expand[c P[x, y]];
       P[y_, c_ * x : (x | a | b)__] := Expand[c P[y, x]];
       P[x_Plus, y_] := P[#, y] & /@ x;
       P[x_, y_Plus] := P[x, #] & /@ y;

In[◦]:= B[y_, x_] := Expand[-B[x, y]];
```

The default value of $n$ (can be changed):

```
In[◦]:= n = 5;
```

The "length" $\lambda$ and the "truth indicator" $\chi_\epsilon$, and the Kronecker $\delta$-function $\delta$:

$$\lambda[x_{i\_,j\_}] := \begin{cases} j - i & i < j \\ n - (i - j) & i > j \end{cases}$$

$$\chi_{\epsilon\_}[cond\_] := \text{If}[\text{TrueQ}@cond, 1, \epsilon];$$

$$\delta_{i\_,j\_} := \chi_0[i == j];$$

The bracket:

$$B[x_{i\_,j\_}, x_{k\_,l\_}] := \begin{cases} \chi_\epsilon[\lambda[x_{i,j}] + \lambda[x_{k,l}] < n]\,(\delta_{j,k} x_{i,l} - \delta_{l,i} x_{k,j}) & j \ne k \vee l \ne i \\ \frac{1}{2} b_i - \frac{1}{2} b_j + \frac{\epsilon}{2} a_i - \frac{\epsilon}{2} a_j & j == k \wedge l == i \end{cases};$$

$$B[a_{i\_}, x_{j\_,k\_}] := (\delta_{i,j} - \delta_{i,k})\, x_{j,k};$$
$$B[b_{i\_}, x_{j\_,k\_}] := \epsilon\,(\delta_{i,j} - \delta_{i,k})\, x_{j,k};$$
$$B[(a | b)\_, (a | b)\_] = 0;$$

The duality pairing:

```mathematica
In[ ]:= P[x_{i_,j_}, x_{k_,l_}] := δ_{j,k} δ_{l,i};
    P[x__, (a | b)_] = 0; P[(a | b)_, x__] = 0;
    P[a_{i_}, b_{j_}] := 2 δ_{i,j}; P[b_{j_}, a_{i_}] := 2 δ_{i,j};
    P[a_, a_] = 0; P[b_, b_] = 0;
```

The permutation $\psi$ and the automorphism $\Psi$:

```mathematica
In[ ]:= ψ[k_Integer] := { k+1   k < n ;
                          1     k == n
    Ψ[ℰ_] := ℰ /. {x_{i_,j_} :> x_{ψ@i,ψ@j}, a_{i_} :> a_{ψ@i}, b_{i_} :> b_{ψ@i}}
```

The basis of $lu_n / gl^{\epsilon}_{n+}$:

```mathematica
In[ ]:= Basis[n_] := Flatten@{
        Table[{x_{i,j}, x_{j,i}}, {i, n-1}, {j, i+1, n}],
        Table[{a_i, b_i}, {i, n}] }
```

## Testing.

```mathematica
In[ ]:= Basis[4]
```

Out[ ]= $\{x_{1,2}, x_{2,1}, x_{1,3}, x_{3,1}, x_{1,4}, x_{4,1}, x_{2,3}, x_{3,2}, x_{2,4}, x_{4,2}, x_{3,4}, x_{4,3}, a_1, b_1, a_2, b_2, a_3, b_3, a_4, b_4\}$

A full bracket-table for $n = 2$:

```mathematica
In[ ]:= n = 2; MatrixForm[
    Table[B[u, v], {u, Basis[n]}, {v, Basis[n]}],
    TableHeadings → {Basis[n], Basis[n]}]
```

Out[ ]//MatrixForm=

|         | $x_{1,2}$ | $x_{2,1}$ | $a_1$ | $b_1$ | $a_2$ | $b_2$ |
|---------|---|---|---|---|---|---|
| $x_{1,2}$ | 0 | $\frac{\epsilon a_1}{2} - \frac{\epsilon a_2}{2} + \frac{b_1}{2} - \frac{b_2}{2}$ | $-x_{1,2}$ | $-\epsilon x_{1,2}$ | $x_{1,2}$ | $\epsilon x_{1,2}$ |
| $x_{2,1}$ | $-\frac{\epsilon a_1}{2} + \frac{\epsilon a_2}{2} - \frac{b_1}{2} + \frac{b_2}{2}$ | 0 | $x_{2,1}$ | $\epsilon x_{2,1}$ | $-x_{2,1}$ | $-\epsilon x_{2,1}$ |
| $a_1$ | $x_{1,2}$ | $-x_{2,1}$ | 0 | 0 | 0 | 0 |
| $b_1$ | $\epsilon x_{1,2}$ | $-\epsilon x_{2,1}$ | 0 | 0 | 0 | 0 |
| $a_2$ | $-x_{1,2}$ | $x_{2,1}$ | 0 | 0 | 0 | 0 |
| $b_2$ | $-\epsilon x_{1,2}$ | $\epsilon x_{2,1}$ | 0 | 0 | 0 | 0 |

The bracket is anti-symmetric at $n = 4$:

```mathematica
In[ ]:= n = 4; Short@Table[
    {u, v} = t; B[u, v] + B[v, u],
    {t, Tuples[Basis[n], 2]}]
```

Out[ ]//Short= $\{0, 0, 0, 0, 0, 0, 0, 0, 0, 0, 0, 0, 0, 0, 0, 0, 0,$
$0, \ll 364\gg, 0, 0, 0, 0, 0, 0, 0, 0, 0, 0, 0, 0, 0, 0, 0, 0\}$

The bracket satisfies the Jacobi identity (strictly, we very Jacobi only for $n = 6$, but three basis elements may involve at most 6 distinct indices, so this is general):

```mathematica
In[ ]:= n = 6; DeleteCases[0]@Table[
    {u, v, w} = t; B[u, B[v, w]] + B[v, B[w, u]] + B[w, B[u, v]],
    {t, Tuples[Basis[n], 3]} ]
```

Out[ ]= { }

The pairing is invariant:

```
In[•]:= n = 6; DeleteCases[0]@Table[
    {u, v, w} = t; P[B[u, v], w] + P[v, B[u, w]],
    {t, Tuples[Basis[n], 3]}]
Out[•]= {}
```

The action of $\Psi$:

```
In[•]:= (# → Ψ[#]) & /@ Basis[4]
```

Out[•]= $\{x_{1,2} \to x_{2,3}, x_{2,1} \to x_{3,2}, x_{1,3} \to x_{2,4}, x_{3,1} \to x_{4,2}, x_{1,4} \to x_{2,1},$
$x_{4,1} \to x_{1,2}, x_{2,3} \to x_{3,4}, x_{3,2} \to x_{4,3}, x_{2,4} \to x_{3,1}, x_{4,2} \to x_{1,3}, x_{3,4} \to x_{4,1},$
$x_{4,3} \to x_{1,4}, a_1 \to a_2, b_1 \to b_2, a_2 \to a_3, b_2 \to b_3, a_3 \to a_4, b_3 \to b_4, a_4 \to a_1, b_4 \to b_1\}$

$\Psi$ is an automorphism:

```
In[•]:= n = 4; DeleteCases[0]@Table[
    {u, v} = t; Ψ[B[u, v]] - B[Ψ[u], Ψ[v]],
    {t, Tuples[Basis[n], 2]}]
Out[•]= {}
```

$\Psi$ respects the pairing:

```
In[•]:= n = 4; DeleteCases[0]@Table[
    {u, v} = t; Ψ[P[u, v]] - P[Ψ[u], Ψ[v]],
    {t, Tuples[Basis[n], 2]}]
Out[•]= {}
```

## Bonus Tests

Acting by arbitrary index-permutations:

```
In[•]:= Act_{σ_List}[ε_] := ε /. {x_{i_,j_} :> x_{σ[[i]],σ[[j]]}, a_{i_} :> a_{σ[[i]]}, b_{i_} :> b_{σ[[i]]}}
```

At $n = 5$, only cyclic permutations induce automorphisms:

```
In[•]:= n = 5;
Select[Permutations[Range[n]],
    σ ↦ And @@ Flatten[Table[
        Act_σ[B[u, v]] === B[Act_σ[u], Act_σ[v]],
        {u, Basis[n]}, {v, Basis[n]}
    ]]
]
Out[•]= {{1, 2, 3, 4, 5}, {2, 3, 4, 5, 1}, {3, 4, 5, 1, 2}, {4, 5, 1, 2, 3}, {5, 1, 2, 3, 4}}
```

Yet in the case of gl$_n$, meaning when $\epsilon = 1$, all permutations induce automorphisms:

```mathematica
In[•]:= n = 4;
    Block[{ϵ = 1}, Select[Permutations[Range[n]],
      σ ↦ And @@ Flatten[Table[
          Act_σ[B[u, v]] === B[Act_σ[u], Act_σ[v]],
          {u, Basis[n]}, {v, Basis[n]}
        ]]
      ]
    ]
Out[•]= {{1, 2, 3, 4}, {1, 2, 4, 3}, {1, 3, 2, 4}, {1, 3, 4, 2}, {1, 4, 2, 3}, {1, 4, 3, 2},
    {2, 1, 3, 4}, {2, 1, 4, 3}, {2, 3, 1, 4}, {2, 3, 4, 1}, {2, 4, 1, 3}, {2, 4, 3, 1},
    {3, 1, 2, 4}, {3, 1, 4, 2}, {3, 2, 1, 4}, {3, 2, 4, 1}, {3, 4, 1, 2}, {3, 4, 2, 1},
    {4, 1, 2, 3}, {4, 1, 3, 2}, {4, 2, 1, 3}, {4, 2, 3, 1}, {4, 3, 1, 2}, {4, 3, 2, 1}}
```

If $\epsilon$ is invertible, the isomorphism class of $\mathfrak{gl}_{n+}^\epsilon$ is independent of $\epsilon$, using Inonu-Wigner contractions:

```mathematica
In[•]:= IW_λ[ℰ_] := ℰ /. {x_{i_,j_} /; i > j ↦ λ x_{i,j}, b_{i_} ↦ λ b_i};

In[•]:= n = 4;
    Union @ Flatten @ Table[
      (B[u, v] /. ϵ → 1) == IW_ϵ @ B[IW_{1/ϵ} @ u, IW_{1/ϵ} @ v],
      {u, Basis[n]}, {v, Basis[n]}]
Out[•]= {True}
```

Even cyclic index permutations become singular at $\epsilon = 0$ when conjugated by Inonu-Wigner contractions (so $\Psi$ simply isn't that):

```mathematica
In[•]:= n = 4; Table[u → IW_{1/ϵ} @ Act_{{2,3,4,1}} @ IW_ϵ @ u, {u, Basis[n]}]
```

Out[•]= $\{x_{1,2} \to x_{2,3},\ x_{2,1} \to x_{3,2},\ x_{1,3} \to x_{2,4},\ x_{3,1} \to x_{4,2},\ x_{1,4} \to \frac{x_{2,1}}{\epsilon},$
  $x_{4,1} \to \epsilon\, x_{1,2},\ x_{2,3} \to x_{3,4},\ x_{3,2} \to x_{4,3},\ x_{2,4} \to \frac{x_{3,1}}{\epsilon},\ x_{4,2} \to \epsilon\, x_{1,3},\ x_{3,4} \to \frac{x_{4,1}}{\epsilon},$
  $x_{4,3} \to \epsilon\, x_{1,4},\ a_1 \to a_2,\ b_1 \to b_2,\ a_2 \to a_3,\ b_2 \to b_3,\ a_3 \to a_4,\ b_3 \to b_4,\ a_4 \to a_1,\ b_4 \to b_1\}$

The same is true for all other permutations (except the identity):

```mathematica
In[•]:= n = 3; MatrixForm[
    Table[IW_{1/ϵ} @ Act_σ @ IW_ϵ @ u, {σ, rows = Permutations @ Range @ n}, {u, cols = Basis[n]}],
    TableHeadings → {rows, cols}]
```

Out[•]//MatrixForm=

|           | $x_{1,2}$ | $x_{2,1}$ | $x_{1,3}$ | $x_{3,1}$ | $x_{2,3}$ | $x_{3,2}$ | $a_1$ | $b_1$ | $a_2$ | $b_2$ | $a_3$ | $b_3$ |
|---|---|---|---|---|---|---|---|---|---|---|---|---|
| {1, 2, 3} | $x_{1,2}$ | $x_{2,1}$ | $x_{1,3}$ | $x_{3,1}$ | $x_{2,3}$ | $x_{3,2}$ | $a_1$ | $b_1$ | $a_2$ | $b_2$ | $a_3$ | $b_3$ |
| {1, 3, 2} | $x_{1,3}$ | $x_{3,1}$ | $x_{1,2}$ | $x_{2,1}$ | $\frac{x_{3,2}}{\epsilon}$ | $\epsilon\, x_{2,3}$ | $a_1$ | $b_1$ | $a_3$ | $b_3$ | $a_2$ | $b_2$ |
| {2, 1, 3} | $\frac{x_{2,1}}{\epsilon}$ | $\epsilon\, x_{1,2}$ | $x_{2,3}$ | $x_{3,2}$ | $x_{1,3}$ | $x_{3,1}$ | $a_2$ | $b_2$ | $a_1$ | $b_1$ | $a_3$ | $b_3$ |
| {2, 3, 1} | $x_{2,3}$ | $x_{3,2}$ | $\frac{x_{2,1}}{\epsilon}$ | $\epsilon\, x_{1,2}$ | $\frac{x_{3,1}}{\epsilon}$ | $\epsilon\, x_{1,3}$ | $a_2$ | $b_2$ | $a_3$ | $b_3$ | $a_1$ | $b_1$ |
| {3, 1, 2} | $\frac{x_{3,1}}{\epsilon}$ | $\epsilon\, x_{1,3}$ | $\frac{x_{3,2}}{\epsilon}$ | $\epsilon\, x_{2,3}$ | $x_{1,2}$ | $x_{2,1}$ | $a_3$ | $b_3$ | $a_1$ | $b_1$ | $a_2$ | $b_2$ |
| {3, 2, 1} | $\frac{x_{3,2}}{\epsilon}$ | $\epsilon\, x_{2,3}$ | $\frac{x_{3,1}}{\epsilon}$ | $\epsilon\, x_{1,3}$ | $\frac{x_{2,1}}{\epsilon}$ | $\epsilon\, x_{1,2}$ | $a_3$ | $b_3$ | $a_2$ | $b_2$ | $a_1$ | $b_1$ |

}

\end{document}